\documentclass[11pt]{amsart}

\usepackage[T1]{fontenc}
\usepackage[utf8]{inputenc}
\usepackage{amsmath,amssymb,amsfonts,amsthm}
\usepackage{mathrsfs}
\usepackage{xspace}
\usepackage{graphicx}
\usepackage{color}
\usepackage{amsmath}
\usepackage{graphicx}
\usepackage{mathrsfs}
\usepackage{amsfonts}
\usepackage{amscd}
\usepackage{wasysym}
\usepackage{color}

\theoremstyle{plain}
\newtheorem{thm}{Theorem}[section]

\newtheorem{cor}[thm]{Corollary}
\newtheorem{con}[thm]{Conjecture}
\theoremstyle{definition}
\newtheorem{defn}[thm]{Definition}
\newtheorem{exam}[thm]{Example}
\theoremstyle{remark}
\newtheorem{rmk}[thm]{Remark}
\newtheorem{Claim}[thm]{Claim}

\newcommand{\g}{\mathfrak g}
\newcommand{\h}{\mathfrak h}
\newcommand{\n}{\mathfrak n}
\newcommand{\bb}{\mathfrak b}
\newcommand{\RR}{\mathbb R}
\newcommand{\CC}{\mathbb C}
\newcommand{\NN}{\mathbb N}
\newcommand{\Cal}[1]{\mathcal{#1}}
\newcommand{\Ad}{\operatorname{Ad}}
\newcommand{\ad}{\operatorname{ad}}
\newcommand{\Ann}{\operatorname{Ann}}
\newcommand{\Ind}{\operatorname{Ind}}
\newcommand{\VA}{V\!A}
\newcommand{\wt}{\widetilde}
\newcommand{\clos}[1]{\overline{#1}^{\,\mathrm{Zar}}}

\begin{document}
	
	\title[The Zariski closure conjecture]{A proof of the Zariski closure conjecture and quasi-Frobenius  exponential algebras}
	
	
	\author{Ali Baklouti}
	\address{Facult\'e des Sciences de Sfax}
	\email{Ali.Baklouti@fss.usf.tn}
	
	\date{October 2026}
	
	\begin{abstract}
		Let $G=\exp(\g)$ be an exponential solvable Lie group, $\pi\in\widehat G$, and $\Omega_\pi$ the coadjoint orbit attached to $\pi$ by the orbit method. The Zariski closure conjecture asserts that the Poisson characteristic variety $\VA(\pi_\nu)$ of the deformed representation coincides with the Zariski closure of $\Omega_\pi$ in $\g^*$. Known until now only in special situations, including the nilpotent case, it is proved here in full generality. Since the Hilbert space of $\pi$ is a quantization of the symplectic manifold $\Omega_\pi$, the statement says that the passage back from the representation to the orbit, a dequantization, loses nothing: the orbit is recovered, up to Zariski closure, from the algebra alone. We also characterize the quasi-Frobenius exponential Lie algebras, those whose dual vector space is a Zariski dense coadjoint orbit, through the triviality of the annihilator of an irreducible representation.
	\end{abstract}
	
	\maketitle
	
	
	
	\section{Introduction}\label{sect:intro}

	The {orbit method} of Kirillov, Bernat and Vergne is one of the most elegant principles of representation theory. It says that the unitary dual $\widehat G$ of a simply connected solvable group should be read off from geometry: irreducible unitary representations correspond to coadjoint orbits in $\g^*$. An orbit $\Omega$ is a symplectic manifold, and the Hilbert space of the corresponding representation is its \emph{quantization}, in the sense of physics. In this dictionary a classical phase space and a quantum state space are two faces of one object.
	
	A dictionary is useful only if it can be read in both directions. The direction from orbits to representations is classical: induce a unitary character from a polarization. The converse, {dequantization}, asks how to recover the orbit from the representation. One concrete answer is algebraic. Put the representation in a family indexed by a deformation parameter $\nu$, let $\nu\to 0$, and look at what survives of the annihilator of the representation in the commutative algebra $S(\g_{\CC})$ of polynomials on $\g^*$. This produces an algebraic subvariety of $\g^*$, the {Poisson characteristic variety} $\VA(\pi_\nu)$. The conjecture states that it is nothing but the Zariski closure of the orbit we started with.   It was formulated in 2005, in joint work with Sami Dhieb and Dominique Manchon \cite{BDM05}, and it is for exponential solvable Lie groups, up-to-date, an open question.\\

	\begin{con}[Zariski closure conjecture]\label{C}
		Let $G$ be an exponential solvable Lie group, and $\pi\in\widehat G$ associated by the orbit method to $\Omega$. Then the Poisson characteristic variety $\VA(\pi_\nu)$ coincides with the Zariski closure of $\Omega$ in $\g^*$.
		
	\end{con}

	The conjecture was born from a several way exchanges; the nilpotent and low-dimensional cases came from patient calculations, and many restrictive solvable settings confirm the dequatization process it means.
	What is at stake mathematically, 
	is that the statement of  Conjecture \ref{C} ties together three worlds that are usually studied separately: the \emph{geometry} of coadjoint orbits, the \emph{analysis} of induced representations, and the \emph{algebra} of primitive ideals in enveloping algebras. In the nilpotent case it follows from the explicit description of annihilators due to Godfrey and Pedersen \cite{Go,Pe84}, where orbits are Zariski closed and $\VA(\pi_\nu)=\Omega$. In the exponential case orbits are no longer closed in the Zariski topology, a spiral being the typical picture (see Example~\ref{Example 1}), and one cannot hope for generators of the annihilator as explicit as in the nilpotent case. The conjecture predicts that, in spite of this, the algebraic object still sees the orbit up to Zariski closure.
	
	Several other partial results were obtained in \cite{BDM19}: the case $\g=[\g,\g]+\g(\ell)$, where $\Omega=G\cdot \ell$;    the case where $\pi$ is induced from of a normal Pukanszky polarization; and a series of low-dimensional examples. An induction on the dimension, through a codimension-one ideal containing the nilradical and $\g(\ell)$, failed: it requires a form of Kirillov's surjectivity theorem which holds in the nilpotent setting and breaks down beyond it. At one point, doubt crept in as well: the successive approximations built term by term in $\nu$ did not seem to converge, and the question of whether a counterexample existed was seriously asked. Failure is part of the story, and it explains why the final argument looks so different from the early attempts: instead of constructing deformations by hand, one uses the generators of the primitive ideals provided by \cite{Bak-Fuji-Casimir} and lets the deformation parameter do the work.
	
	The main result of this note is a proof of the above conjecture for an arbitrary exponential solvable Lie group.  We have the following:

	\begin{thm}\label{main result}
		Conjecture \ref{C} holds for an exponential solvable Lie group $G$.  For any  $\pi\in\widehat G$ with associated coadjoint orbit $\Omega=\Omega(\pi)$, 
		$\VA(\pi_\nu)=\clos{\Omega}$, the Zariski closure of $\Omega$ in $\g^*$.
	\end{thm}
	
	The coadjoint orbits of an exponential Lie group are not algebraic  smooth  varieties as in the nilpotent setting. However, the proof is short as the work is done upstream, and it is given in Section~\ref{sect:proof}.  It  is based on a theorem of \cite{Bak-Fuji-Casimir} which provides a family of generators of $\Ann\pi$ whose symbols generate the ideal of polynomials vanishing on $\Omega$, and the compatibility of these generators with the deformation $\g\mapsto\g_\nu$ does the rest. 
	
	What remains, once the technicalities are set aside, is the image from the beginning: a representation and an orbit are two descriptions of the same object, and Theorem \ref{main result} says that no information is lost, up to a natural closure, in going from one to the other. That such a statement could be reached only after many years, through the work of several people on different questions, and consequent  shared endeavour may lead to very important conclusions.

	In Section~\ref{sect:QF} we single out the  class of  quasi-Frobenius Lie algebras, and provide a complete characterization in the exponential setting by means of primitive ideals.   Section~\ref{sect:persp} is devoted to indicate open directions.
	
	\medskip

	\section{The setting: orbits, representations, deformation}\label{sect:setting}
	
	We recall the facts and notation needed later; see \cite{LL94} for details on exponential groups.
	
	\subsection{Exponential groups and the orbit method}
	
	Let $G$ be a connected real Lie group with Lie algebra $\g$. It is \emph{exponential solvable} if $\exp:\g\to G$ is a $C^\infty$-diffeomorphism; in particular $G$ is simply connected. The group acts on $\g$ by $\Ad$ and on the dual $\g^*$ by
	\[
	\Ad^*(g)\ell(Y)=\ell(\Ad(g^{-1})Y),\qquad g\in G,\ \ell\in\g^*,\ Y\in\g .
	\]
	The coadjoint orbit of $\ell$ is $\Omega_\ell=\{g\cdot\ell:\ g\in G\}$, and $\g^*/G$ denotes the orbit space.
	
	For $\ell\in\g^*$ let $$B_\ell(X,Y)=\ell([X,Y])$$ with radical $\g(\ell)$. A \emph{polarization} at $\ell$ is a subalgebra $\bb\subset\g$, totally isotropic for $B_\ell$ and of maximal dimension $\tfrac12(\dim\g+\dim\g(\ell))$. It satisfies \emph{Pukanszky's condition} if
	\[
	\Ad^*(\exp\bb)\,\ell=\ell+\bb^\perp .
	\]
	Put $B=\exp\bb$, let $\chi_\ell(\exp X)=e^{i\ell(X)}$ for $X\in\bb$ and
	\[
	\pi_{\ell,\bb}=\Ind_B^G\chi_\ell .
	\]
	The Hilbert space of $\pi_{\ell,\bb}$ is realized as a space of functions on $G/B\simeq\RR^{d/2}$, where $d=\dim\Omega_\ell$. This is the sense in which it is a quantization of $\Omega_\ell$. The theorem of Kirillov--Bernat--Vergne states that $\pi_{\ell,\bb}$ is irreducible, that its class depends only on $\Omega_\ell$, and that every $\pi\in\widehat G$ arises in this way; the resulting map
	\[
	\g^*/G\longrightarrow\widehat G,\qquad \Omega_\ell\longmapsto [\pi_{\ell,\bb}],
	\]
	is a homeomorphism for the quotient and Fell topologies \cite{LL94}. For $\pi\in\widehat G$ we write $\Omega(\pi)$ for the orbit corresponding to $\pi$.
	
	\subsection{The deformation}\label{subsect:deformation}
	
	Let $\g_{\CC}$ be the complexification of $\g$ and $S(\g_{\CC})$ its symmetric algebra, that is, the polynomial functions on $\g^*$. For $\nu\in\CC$ let $\g_\nu$ be the vector space $\g$ with bracket
	\[
	[X,Y]_\nu=\nu\,[X,Y].
	\]
	For $\nu\neq0$ this Lie algebra is isomorphic to $\g$; for $\nu=0$ it is abelian. The symmetrization map
	\[
	\sigma_\nu:S(\g)\xrightarrow{\ \sim\ }\mathcal U(\g_\nu),\qquad
	x_1\cdots x_n\longmapsto\frac1{n!}\sum_{\tau\in S_n}x_{\tau(1)}\cdots x_{\tau(n)},
	\]
	is a linear isomorphism and transports the product of $\mathcal U(\g_\nu)$ to a (star)-product on $S(\g)$:
	\begin{equation}\label{star}
		P*_\nu Q:=\sigma_\nu^{-1}\big(\sigma_\nu(P)\,\sigma_\nu(Q)\big).
	\end{equation}
	For all $P,Q$, the expression $P*_\nu Q$ is polynomial in $\nu$, and
	\[
	P*_\nu Q=PQ+O(\nu),\qquad P*_\nu Q-Q*_\nu P=\nu\{P,Q\}+O(\nu^2),
	\]
	where $$\{f,g\}(\ell)=\langle\ell,[df(\ell),dg(\ell)]\rangle, f,g\in S(\g_{\CC}), \ell\in \g^*$$
	is the Kirillov--Kostant--Souriau Poisson bracket. Thus $*_\nu$ is a (Lugo--Gutt type) deformation quantization of the Poisson algebra $S(\g)$. The same construction may be carried out with the Duflo isomorphism $\tau=\sigma\circ J(D)^{1/2}$ in place of $\sigma$ \cite{18}, where $$J(x)^{1/2}=\big(\det\frac{\operatorname{sh}\ad(\nu x/2)}{\ad(\nu x/2)}\big)^{1/2};$$ only the symmetrization will be used below.
	
	\subsection{The module $\pi_\nu$ and the Poisson characteristic variety}
	
	Let $\rho$ be an irreducible unitary representation of $G=\exp\g$ attached to $\Omega\ni\ell$ with Pukanszky polarization $\h$. For every real $\hbar$ the same construction gives a representation $\rho_\hbar$ of $G_\hbar=\exp\g_\hbar$, irreducible for $\hbar\neq0$. Differentiating and setting $\nu=i\hbar$, with $\nu$ regarded as an indeterminate, one obtains the topological module $(\pi_\nu,\Cal M)$ over the formal enveloping algebra, $\pi_\nu=i\rho_{-i\nu}$ \cite[\S3]{BDM05}. Let $A=S(\g_{\CC})$; the module induces an $A$-module structure on $M=\Cal M/\nu\Cal M$, with representation $\pi_0$. Two ideals of $A$ enter:
	\[
	I=\Ann\pi_0,\qquad J=\Ann\pi_\nu/\nu\Ann\pi_\nu .
	\]
	The \emph{characteristic variety} $V(\pi_\nu)$ and the \emph{Poisson characteristic variety} $\VA(\pi_\nu)$ are the common zero sets of $I$ and $J$, respectively, viewed as real varieties by unitarity. Since $J\subset I$, we have $V(\pi_\nu)\subset\VA(\pi_\nu)$. It was proved in \cite[Th\'eor\`eme~5.3.1]{BDM05}, for every exponential solvable group, that
	\[
	V(\pi_\nu)=\ell+\h^\perp .
	\]

	\subsection{Layering and generators of primitive ideals}
	
	Let
	\[
	\mathcal S:\ \g_1\subset\cdots\subset\g_{m-1}\subset\g_m=\g
	\]
	be a composition series of ideals of $\g$ passing through the nilradical $\n=\g_{j_0}$, so that $1\le\dim(\g_j/\g_{j-1})\le2$. For $\ell\in\g^*$ put $e_j(\ell)=\dim G\cdot(\ell|_{\g_j})$ and $e(\ell)=(e_1(\ell),\dots,e_m(\ell))$. The set $\mathcal E=\{e(\ell)\}$ is ordered lexicographically; it defines the layers $U_e=\{\ell: e(\ell)=e\}$, which are $G$-invariant, and the maximal layer is Zariski open. The index sets $S(e)=\{j:\ e_j\neq e_{j-1}\}$ and $T(e)=\{j:\ e_j=e_{j-1}\}$ are the jump and non-jump indices.\\

		Let $e \in \mathcal E$. We say that $A \in \mathcal U(\mathfrak g)$ is an
\textit{	$e$-Casimir}
	element (cf.  \cite{Bak-Fuji-Casimir}), if there exists a non-empty $G$-invariant Zariski open set $\mathcal Z$ of ${\mathfrak g}^*$
	such that $\Sigma_e = \mathcal Z \cap U_e$ is non-empty, and that ${\pi}_{\ell}(A)$ is a scalar operator for $\ell \in \Sigma_e $.   We note $\operatorname{Cas}(\Sigma_e )$ the algebra of $e$-Casimir elements, and we write ${\pi}_{\ell}(A)
	= {\varphi}_A(\ell)I$ and we denote {by} ${\varphi}_A(\ell)$ the $e$-Casimir function on $\Sigma_e$ corresponding to the $e$-Casimir element $A$.
	
	We extract from to $\mathcal S$,  a good sequence of subalgebras $\mathcal S':\ \mathfrak k_1\subset\dots\subset\mathfrak k_n=\g$, which means that $\dim\mathfrak k_i=i$ and if $\mathfrak k_i$ is not an ideal of $\mathfrak g$, then $\mathfrak k_{i-1}$ and $\mathfrak k_{i+1}$ are.  Let  $$\bb[\ell]=\sum_{i=1}^{n}\mathfrak k_i(\ell\vert_{\mathfrak k_i})$$ be the associated Vergne polarization, which allows to define $\pi_\ell=\Ind_{B[\ell]}^G\chi_\ell$, where $B[\ell]=\exp(\bb[\ell])$. Let $\mathcal K$ be the set of indices $j_0<k\le m$ with $e_k=e_{k-1}+1$, so that $\bb[\ell]\cap\g_k=\bb[\ell]\cap\g_{k-1}$ for generic $\ell$, and write $\g_k=\RR X_k+\g_{k-1}$ for $k\in\mathcal K$.
	
	We now set the following:
	
	\begin{defn}[cf.~\cite{Bak-Fuji-Casimir}]\label{def:J}
		Let $\mathcal J$ be the set of indices $1\le j\le m-1$ such that
		\begin{enumerate}
			\item $e_j\neq e_{j-1}+\dim(\g_j/\g_{j-1})$;
			\item there exist $W\in\mathcal U(\g_j)\setminus\mathcal U(\g_{j-1})$ in $\Ann\pi|_{G_{m-1}}$ and some integers $n_k\in\NN$ with $(\ad X_k)^{n_k}(W)\in\Ann\pi$ for all $k\in\mathcal K$.
		\end{enumerate}
	\end{defn}

	We then have the following result:

	\begin{thm}[{\cite[Thm.~4.6, Prop.~4.7]{Bak-Fuji-Casimir}}]\label{generators of the kernel}
		Let $G=\exp\g$ be an exponential solvable Lie group, $\pi\in\widehat G$ and $\Omega=\Omega(\pi)$. With the notation above:
		\begin{enumerate}
			\item in passing from $\g_{j-1}$ to $\g_j$ there is a new element $W_j\in\mathcal U(\g_j)\cap\Ann\pi$ if and only if $j\in\mathcal J$;
			\item the primitive ideal $\Ann\pi$ is generated by $(W_j)_{j\in\mathcal J}$;
			\item the Poisson ideal $I(\Omega)$ of polynomials on $\g^*$ vanishing on $\Omega$ is generated by $Q_j=\sigma^{-1}(W_j)$, $j\in\mathcal J$.
		\end{enumerate}
	\end{thm}
	
	This result is the engine of the proof of Theorem \ref{main result}. 
	The mechanism behind the construction of the generators of $\Ann\pi$ proceeds by induction on $\dim\g$, following the composition series of ideals $(\mathcal S)$. At each passage from $\g_{j-1}$ to $\g_j$ one asks whether $\Ann\pi\cap\mathcal U(\g_j)$ contains a \emph{new} element not already in $\mathcal U(\g_{j-1})$. The answer is governed by two conditions, which define the index set $\mathcal J$. The first is $e_j\neq e_{j-1}+\dim(\g_j/\g_{j-1})$. When equality holds, the orbit of $\ell|_{\g_j}$ is saturated with respect to $\g_{j-1}$, the restriction of $\pi$ is a direct integral $\int^\oplus\rho_t\,dt$ of pairwise inequivalent representations over the transversal variable, and no new relation can appear at this level. The second condition asks for $W\in\mathcal U(\g_j)\setminus\mathcal U(\g_{j-1})$ lying in $\Ann\pi|_{G_{m-1}}$ and such that $(\operatorname{ad}X_k)^{n_k}W\in\Ann\pi$ for the indices $k$ with $e_k=e_{k-1}+1$. This second condition expresses that the relation survives the transport along the orbit directions.
	
	The proof splits according to the position of $\Omega(\pi)$ with respect to the codimension-one ideal $\g_{m-1}$. If $\Omega$ is \emph{saturated}, then $\pi\simeq\operatorname{Ind}_{G_{m-1}}^G\rho_t$ and $\pi|_{G_{m-1}}\simeq\int^\oplus\rho_t\,dt$. Since $\pi(X_m)$ acts as a derivation in the transversal variable, any $W=\sum_k a_kX_m^k\in\Ann\pi$ has all coefficients $a_k$ in $\Ann\rho_t$, so the induction hypothesis applies to $\g_{m-1}$ and the kernel is generated by the family  of generators attached to $\rho_t$.

	If $\Omega$ is \emph{non-saturated}, one passes to the quotient by any ideal contained in $\Ann\pi$, so that the center $\mathfrak z $ of $\mathfrak g$ verifies $\dim\mathfrak z\le1$, and one takes a minimal non-central ideal $\mathfrak a\supset\mathfrak z$ (abelian of dimension $\le3$). In this case, Proposition~4.2 in  \cite{Bak-Fuji-Casimir}, provides $W_0=X_m+Q\in\Ann\pi$ with $Q\in\mathcal U(\g_{m-1})$. A relation of degree one in the extra variable makes $\pi|_{G_{m-1}}$ irreducible, the remaining generators live in $\mathcal U(\mathfrak a^\ell)$ with $\mathfrak a^\ell=\{T\in \frak g:\ \ell([T,\g])=0\}\neq\g$, and the induction applies again to produce the family of generators $\mathscr A$ of $\Ann \pi$.\\

	The passage to the polynomial ideal rests on one observation. For a Casimir element, $\pi_\ell(W)=\varphi_W(\ell)\,I=\sigma^{-1}(W)(\ell)\,I$ (cf. \cite{Bak-Fuji-Casimir}, Proposition~3.1), so $A_j=W_j-\varphi_{W_j}(\ell)$ belongs to $\Ann\pi$ and its symbol $\sigma^{-1}(A_j)$ vanishes on $\Omega(\pi)$. In the saturated case the same conclusion follows from the direct-integral decomposition and the induction hypothesis, and in the non-saturated case it follows from the explicit form of the elements in lower stages.   Consequently $\sigma^{-1}(\mathscr A)$ generates $I(\Omega(\pi))$.  
	
	Finally, Theorem~4.9 in \cite{Bak-Fuji-Casimir} shows that the generators of $\Ann \pi$ may be chosen to be Casimir elements: $W=u-\varphi_u$ ($\varphi_u$ denotes the $e$-Casimir function)  with $u\in\operatorname{Cas}(\Sigma_e)$, since a direct computation on $L^2(\mathbb R,\mathcal H_{\pi'})$ gives $\pi_\ell(u)=\varphi_u(\ell')\,I$. Thus the passage from representation to orbit is carried by a family of scalar-valued operators, one for each new index $j\in\mathcal J$.

	\section{Proof of the conjecture}\label{sect:proof}
	
	\subsection{Reformulation}
	
	Recall the descriptions of the two ideals from \cite{BDM05}: $I$ consists of the $Q_0\in S(\g_\CC)$ such that $\pi_\nu(Q_0)=O(\nu)$, whereas $J$ consists of the $Q_0$ for which there exist polynomials $Q_k$, $k\ge1$, with
	\[
	Q=Q_0+\nu Q_1+\nu^2Q_2+\cdots\in\Ann\pi_\nu .
	\]
	In other words, $I$ is the set of polynomials whose image is ``small'' at $\nu=0$, and $J$ the set of polynomials which can be \emph{deformed into the annihilator}. One has $J\subset I$.
	
	\begin{rmk}[The easy inclusion]\label{easy}
		$\clos{\Omega}\subset\VA(\pi_\nu)$ always holds. Indeed, if $Q_0\in J$ is deformed into $Q\in\Ann\pi_\nu$, then for $X\in\g$ the commutator $X*_\nu Q-Q*_\nu X=\nu\{X,Q_0\}+O(\nu^2)$ lies in the two-sided ideal $\Ann\pi_\nu$; dividing by $\nu$ (the module is torsion free) shows $\{X,Q_0\}\in J$. Hence $J$ is stable under the coadjoint action, so is its zero set $\VA(\pi_\nu)$, which contains $V(\pi_\nu)=\ell+\h^\perp$; by Pukanszky's condition, the $G$-saturation of this set is $\Omega$. Being closed, $\VA(\pi_\nu)$ then contains $\clos\Omega$.
	\end{rmk}
	
	Consequently Conjecture~\ref{C} is equivalent to the reverse inclusion, which amounts to the following claim.
	\begin{Claim}\label{Cbis}
		Any polynomial vanishing on the orbit can be deformed into an element of the annihilator of the topological module $\pi_\nu$; that is, $I(\Omega)\subset J$.
	\end{Claim}
	
	\subsection{The deformation argument}
	
	\begin{proof}[Proof of Theorem~\ref{main result}]

		Applying  Theorem \ref{generators of the kernel} to $\g_\nu$ for $\nu\neq0$ (an exponential solvable Lie algebra, isomorphic to $\g$), one obtains generators $W_{j,\nu}$ of $\Ann\pi_\nu$ and, applying $\sigma_\nu^{-1}$, generators
		\[
		Q_{j,\nu}=\sigma_\nu^{-1}(W_{j,\nu}),\qquad j\in\mathcal J,
		\]
		of the Poisson ideal of $\big(S(\g),*_\nu\big)$. For $\nu\ne0$ the orbit of $\ell$ under the group of $\g_\nu$ is again $\Omega$, so this ideal is $I(\Omega)$, and the $Q_{j,\nu}$ are \emph{polynomial in $\nu$}. We write ${\widetilde Q_j}:=Q_{j,0}$; thus
		\begin{equation}\label{Qnu}
			Q_{j,\nu}={\widetilde Q_j}+O(\nu),\qquad \sigma_\nu(Q_{j,\nu})\in\Ann\pi_\nu,\qquad I(\Omega)=\big\langle {\widetilde Q_j}:\ j\in\mathcal J\big\rangle.
		\end{equation}
		Let $P\in I(\Omega)\subset S(\g)$. In the form recorded in \eqref{Qnu}, there are $R_j\in S(\g)$, $j\in\mathcal J$, with
		\[
		P=\sum_{j\in\mathcal J}{\widetilde Q_j}R_j .
		\]
		Consider
		\[
		\wt P_\nu:=\sum_{j\in\mathcal J}\sigma_\nu(Q_{j,\nu})\,\sigma_\nu(R_j)\ \in\ \mathcal U(\g_\nu).
		\]
		It is polynomial in $\nu$, and it belongs to the two-sided ideal $\Ann\pi_\nu$, since each $\sigma_\nu(Q_{j,\nu})$ does. Put $P_\nu:=\sigma_\nu^{-1}(\wt P_\nu)=\sum_jQ_{j,\nu}*_\nu R_j$. Because $*_\nu$ is the ordinary product modulo $\nu$ and $Q_{j,\nu}={\widetilde Q_j}+O(\nu)$, we get
		\[
		P_\nu=\sum_{j\in\mathcal J}{\widetilde Q_j}R_j+O(\nu)=P+O(\nu),
		\]
		so $P_\nu$ is a deformation of $P$ lying in $\Ann\pi_\nu$. Hence $P\in J$. This proves $I(\Omega)\subset J$, that is $\VA(\pi_\nu)=V(J)\subset V(I(\Omega))=\clos\Omega$. Together with Remark~\ref{easy}, $\VA(\pi_\nu)=\clos\Omega$.
	\end{proof}
	
	\begin{rmk}
		In the nilpotent case orbits are Zariski closed and the theorem gives $\VA(\pi_\nu)=\Omega$, recovering \cite[Th\'eor\`eme~5.4.1]{BDM05}. There the generators of $J$ were also generators of $\Ann\rho$ after symmetrization, via the descriptions of Godfrey and Pedersen \cite{Go,Pe84}. The argument above shows that this feature, which looked specific to nilpotent groups, is shared by all exponential solvable groups.
	\end{rmk}
	
	\section{Quasi-Frobenius  Lie algebras}\label{sect:QF}
	
	A Lie algebra $\g$ is \emph{Frobenius} if some $\ell\in\g^*$ makes $B_\ell$ non-degenerate, equivalently if $\g$ admits an open coadjoint orbit. The notion goes back to Ooms \cite{O74,O80}, who studied it through primitive ideals: the zero ideal of a Frobenius algebra is primitive, but not conversely. For a general Lie algebra, an orbit $\Omega=G\cdot\ell$ is open in an affine subspace if and only if $\g(\ell)$ is an ideal of $\g$ \cite[Prop.~3.1]{Bak-Ishi}. For exponential solvable algebras one has the following characterization \cite[Thm.~1.1]{Bak-Ishi}.
	
	\begin{thm}\label{char Frob Exp}
		For an exponential solvable Lie algebra $\g$, the following are equivalent:
		\begin{enumerate}
			\item $\g$ is a Frobenius Lie algebra;
			\item there is $\ell\in\g^*$ such that $\g(\ell)$ is an ideal of $\g$ and $\Ann \pi_\ell=\{0\}$;
			\item there is $\ell\in\g^*$ such that $\Omega_\ell$ is relatively open in an affine subspace of $\g^*$ and $\Ann\pi_\ell=\{0\}$.
		\end{enumerate}
	\end{thm}
	
	In particular, if $\g$ is Frobenius, some $\pi_\ell\in\widehat G$ has $\Ann \pi_\ell=\{0\}$. The converse fails, as the following example shows.
	
	\begin{exam}\label{Example 1}
		Let $\g_3$ be spanned by $A,X,Y$ with $[A,X]=X-Y$, $[A,Y]=X+Y$, and $G=\exp\g_3$, an exponential, non completely solvable group. Let $f=xX^*+yY^*+aA^*\in\g_3^*$. If $x^2+y^2=0$ the orbit of $f$ is a point. Otherwise $\bb=\RR X\oplus\RR Y$ is a polarization satisfying Pukanszky's condition; let $\chi_f(\exp U)=e^{-if(U)}$ on $B=\exp\bb$ and $\rho_f=\Ind_B^G\chi_f$. Then there is a unique $\theta\in[0,2\pi[$ with $\rho=\rho_\theta=\rho_{f_\theta}$, $f_\theta=\cos\theta\,X^*+\sin\theta\,Y^*$, and its orbit is
		\[
		\Omega_\theta=\big\{sA^*+e^{-t}\cos(t+\theta)X^*+e^{-t}\sin(t+\theta)Y^*:\ s,t\in\RR\big\}.
		\]
		Moreover,
		\[
		d\rho(A)=-\frac{d}{dt},\qquad d\rho(X)=-i\,e^{t}\cos(\theta+t),\qquad d\rho(Y)=-i\,e^{t}\sin(\theta+t),
		\]
		and one checks that the annihilator of $\rho$ is trivial. Yet the orbit $\Omega_\theta$ is not open (it has dimension $2$ in a $3$-dimensional space): its projection to the $(X^*,Y^*)$-plane is a logarithmic spiral, which is Zariski dense in the plane. Thus $\Omega_\theta$ is Zariski dense in $\g_3^*$.
	\end{exam}
	
	This example is automatically non-Frobenius because of its odd dimension; a $4$-dimensional counterexample is given in \cite[\S4]{Bak-Ishi}. It leads us to isolate the right notion.
	
	\begin{defn}\label{Quasi-Frobenius}
		A real Lie algebra $\g$ is \emph{quasi-Frobenius} if $\g^*$ is a quasi-orbit, that is, if there is a Zariski dense coadjoint orbit.
	\end{defn}
	
	Any Frobenius Lie algebra is quasi-Frobenius,  and the algebra $\g_3$ above is quasi-Frobenius but not Frobenius, while no nilpotent Lie algebra is.
	
	\begin{cor}\label{quasi-Frobenius characterization}
		An exponential Lie algebra $\g$ is quasi-Frobenius if and only if $G=\exp\g$ admits $\pi_\ell\in\widehat G$ with $\Ann\pi_\ell=\{0\}$ for some $\ell\in\g^*$.
	\end{cor}
	
	\begin{proof}
		By Theorem~\ref{generators of the kernel}, $\Ann\pi_\ell$ is generated by elements $W_j$ whose symbols $Q_j=\sigma^{-1}(W_j)$ generate $I(\Omega_\ell)$. Hence $\Ann\pi_\ell=\{0\}$ if and only if $I(\Omega_\ell)=\{0\}$, which means that no non-zero polynomial vanishes on $\Omega_\ell$, that is, $\Omega_\ell$ is Zariski dense in $\g^*$.
	\end{proof}
	
	In the language of the main theorem: $\g$ is quasi-Frobenius precisely when some $\pi$ has $\VA(\pi_\nu)=\g^*$. The Frobenius algebras, with their open orbits, are then a special and rigid subclass, while the quasi-Frobenius ones, which include the spiral of Example~\ref{Example 1}, form the natural class on the ``non-algebraic'' side of the theory.
	
	\section{Perspectives}\label{sect:persp}

	\subsection{Setting of complex solvable Lie algebras}

	Let now $\frak s$ be a complex solvable Lie algebra,
	and $\mathrm{Prim}( \mathcal U(\frak s))$
	the set of primitive ideals of $\mathcal U (\frak s)$.
	Let $\mathcal G$  be the algebraic adjoint of $\frak s$.
	Namely, 
	$\mathcal G$ is the smallest algebraic subgroup of $GL(\frak s)$
	whose Lie algebra contains $\mathrm{ad}(\frak s) \subset \mathfrak{gl}(\frak s)$.
	Then the  group $\mathcal G$ acts on $\frak s^*$ by the contragredient representation.
	Dixmier \cite{D66} introduced the  map:
	\begin{eqnarray*}
		\mathcal D:\mathfrak s^*/\mathcal G&\longrightarrow& \rm{Prim}\ \mathcal U(\mathfrak s)\\
		\hbox{Ad}^*\mathcal G.\ell &\longmapsto & \mathcal D(\ell),
	\end{eqnarray*}
	which associates to any $\mathcal G$-orbit $\hbox{Ad}^*\mathcal G.\ell$ in $\mathfrak s^*$ a primitive ideal of the enveloping algebra $\mathcal U(\mathfrak s)$. 
	This map is invariant under this action,
	which means that 
	$\mathcal D(\hbox{Ad}^*g .\ell) = \mathcal D(\ell)$ 
	for $g \in \mathcal G$ and $\ell \in \frak s^*$.

	The Dixmier map $\mathcal D$ has been shown to be injective by M. Duflo \cite{Du72}, and surjective by R. Rentschler \cite{R74}. Continuity of the map $\ell\mapsto \mathcal D(\ell)$ was proven by N. Conze and M. Duflo \cite{CD70} as early as 1970 (when $\mathfrak s^*$ is endowed with the Zariski topology, and $\mathrm{Prim}\ \mathcal U(\mathfrak s)$ with the Jacobson topology). Finally, the bi-continuity of the Dixmier map has been established by O. Mathieu in 1991 (cf. \cite{M91}).

	As an insight on the relation
	between primitive ideals and algebraic adjoint orbits has also been investigated in \cite{Bak-Ishi}.  
	We showed the following:
	%
	%
	
	\vskip 0.4cm
	
	\begin{thm} \label{thm:OpenO}
		{For a complex solvable Lie algebra $\frak s$ and $\ell \in \frak s^*$ with the notation above},
		the primitive ideal $\mathcal D(\ell)$ is zero if and only if
		the orbit $\hbox{Ad}^*\mathcal G.\ell  \subset \frak s^*$ is open.
	\end{thm} 
	Denoting by $J_\ell\subset S(\mathfrak s)$ the (Poisson) ideal of polynomials vanishing on $\rm{Ad}^*\mathcal G.\ell$, a key question is to investigate  whether it is possible  to deform any $P\in J_\ell$ to get an element $\wt P\in \mathcal D(\ell)$.

	\subsection{Characteristic varieties in the real solvable case}  One further question seems natural.  Is it possible to go beyond the exponential class to more general solvable groups, where the situation is considerably more subtle? In this setting, real polarizations do not suffice to construct irreducible representations, and on can assign to a given coadjoint orbit an infinite  number of distinct inequivalent
	representations. So the orbit method 
	itself has to be modified.

	\subsection{$C^*$-algebraic version of characteristic varieties}
	
	The Kirillov--Bernat map $\g^*/G\to\widehat G$ is a homeomorphism by Leptin and Ludwig \cite{LL94}, and $\widehat G$ is identified with $\operatorname{Prim}C^*(G)$ since exponential groups are of type~I (see \cite{BDM11}). Our theorem says that, after Zariski closure, the inverse of the Kirillov map is computed by dequantization: choose $\rho$, differentiate the family $\rho_\hbar$, and read off $\VA$. It would be interesting to study the bicontinuity of the Kirillov map in this light, and to find a $C^*$-algebraic version of $\VA$ based on the tangent groupoid $\coprod_{\hbar\in[0,1]}G_\hbar$, where the closure in the ordinary topology of $\g^*$ should replace the Zariski closure.

	\bigskip
	\noindent\textbf{Acknowledgements.} The author deeply thanks Dominique Manchon, for many discussions about the subject.
	

\end{document}